\documentclass[12pt]{article}

\usepackage{amsmath,amssymb}

\DeclareMathOperator{\tr}{tr}
\DeclareMathOperator{\id}{id}
\DeclareMathOperator{\sh}{sh}

\newcommand{\AAA}{\mathcal{A}}
\newcommand{\BB}{\mathcal{B}}
\newcommand{\CC}{\mathbb{C}}
\newcommand{\CCC}{\mathcal{C}}
\newcommand{\DD}{\mathcal{D}}
\newcommand{\RR}{\mathbb{R}}
\newcommand{\RRR}{\mathcal{R}}
\newcommand{\UU}{\mathcal{U}}
\newcommand{\ZZ}{\mathbb{Z}}

\newcommand{\ONE}{\mbox{I}}
\newcommand{\TWO}{\mbox{II}}

\begin{document}

\fnsymbol{footnote}

\begin{center}
{
\bf \large \bf Modular Double of Quantum Group \\
}

\vskip 2.0cm

{\bf Ludvig Faddeev } \\ 
\vskip 0.3cm

{\it St.Petersburg Branch of Steklov Mathematical
Institute \\
Russian Academy  of Sciences, Fontanka 27 , St.Petersburg,
Russia } \\

\vskip 0.3cm

{\it Helsinki Institute of Physics \\
P.O. Box 9, FIN-00014 University of Helsinki, Finland} \\
\end{center}

\vskip 0.3cm

	As it is clear from the title, I shall deal with some question 
	connected with the theory of Quantum Groups. If I remember right,
	Moshe did not like Quantum Groups after this notion was cristallized
	by Drinfeld 
\cite{Drinfeld}
	in pure algebraic manner. However his own attraction to the
	deformations (as well as the pressure of authors in LMP) made him
	to change his mind. So when I presented the subject described
	below at St. Petersburg meeting on May 1998, he did not express any bad
	feelings. So I decided to publish it in this memorial volume.
	
	There are several sources of my proposal. I shall give just two,
	one "mathematical" and another "physical", as it is appropriate for
	a paper on Mathematical Physics.
	
	1. In the definition of Quantum Group one uses the deformation of
	the Chevalley generators
    $ K $, $ f $, $ e $,	
	whereas for the construction of the universal	
	$ R $-matrix one needs nonpolinomial elements like
    $ H = \ln K $.	
	Explicite formulas will be reminded below. This unfortunate
	obstacle can be circumvented in several ways: one, \`{a}-la Lusztig
\cite{Lusztig}	
	is just not to use explicite formula of Drinfeld; another,
	followed in the most of texts on Quantum Groups (see i.e.
\cite{Kassel}), 
	is to employ formal series in
    $ \ln q $.
	However the value of
    $ R $-matrix as a genuine operator is too high and deserves more friendly
	attitude.
	
	2. In the applications of Quantum Groups to Conformal Field Theory
	one explicitely sees, that together with the attributes of Quantum
	Groups (i.e. eigenvalues of Laplacians) for
    $ q = e^{i \pi \tau} $	
	there enter analogous objects, corresponding to
    $ \tilde{q} = e^{- i \pi / \tau} $.		
	This modular duality is a well known "experimental" fact, which
	goes without satisfactory explanation.
	
	In the following an extension of Quantum Group will be described,
	which will throw some light on both topics above. Roughly speaking
	I propose to unite the Quantum Groups for 
    $ q $ and $ \tilde{q} $	
	in one object having modular structure (see e.g.
\cite{Haag}). 
	Thus the combination of words "Modular double of Quantum Group"
	seems to be quite relevant for it.
	
	The element
    $ \ln K $
	will have a natural definition and the expression for the universal
    $ R $-matrix 
	will become much more meaningfull in this extension. I
	believe also, that it is this modular double which in fact defines
	the hidden symmetry in the Conformal Field Theory. There are some
	indication of this in the literature
\cite{Gervais,Babelon,BLZ}
	and recently it was made more explicite in
\cite{FV,FVK}.
	At this moment I know well how my proposal works in all detail in
	the case of rank 1
    $ SL(2) $
	group. This will be presented below. The tools for the
    $ SL(N) $
	generalizations are known
\cite{Tarasov,KV},
	but other serieses of simple groups are not treated yet.
	
	The work was partially supported by INTAS, 
	RFFR (grant RFFR-96-01-00851)
	and University of Helsinki.
	
\section{Reminder on Quantum Group 
    $ SL(2)$ }
	I shall use an extension of
    $ SL_{q}(2) $	
	which is Drienfeld double of its Borel part. There are four
	generators
    $ K $, $ K' $, $ e $, $ f $		
	with familiar relations
\begin{eqnarray}
    && Ke = q^{2} eK \, ; \quad K'e = q^{-2} eK' \,; \\ 
    && Kf = q^{-2} fK \, ; \quad K'f = q^{2} fK' \,; \\ 
    && ef - fe = \frac{K-K'}{q-q^{-1}} \, ; \quad K K' = K' K \, .
\end{eqnarray}	
	The algebra
    $ \UU_{q} $
	generated by
    $ K $, $ K' $, $ e $, $ f $			
	over the field
    $ \CC $
	of complex numbers, (so that 
    $ q $	
	is a complex number) has two central elements
\begin{equation}
    J = KK' \, ; \quad C = \frac{K-K'}{q-q^{-1}} +
	(q-q^{-1})^{2} (ef-fe) \: .
\end{equation}
	Reduction to 
    $ SL_{q}(2) $		
	is achieved if we put
    $ J = \id $,
	however, we shall not do it here.
	
	The universal
    $R$-matrix is affiliated with the tensor square
    $ \UU_{q} \otimes \UU_{q} $	
	and defines the property of the Hopf multiplication 
    $ \bigtriangleup $
	in
    $ \UU_{q} $
\begin{equation}
    \sigma \circ \bigtriangleup = R \bigtriangleup R^{-1} \, ,
\end{equation}	
	where
    $ \sigma $
	is a permutation
    $ \sigma (a \otimes b) = b \otimes a $.
	Drinfeld has given a formal expression for
    $ R $
\begin{equation}
    R = q^{- \frac{H\otimes H'}{2}} s_{q} (-(q-q^{-1})^{2} e \otimes f )
\end{equation}	
	where
\begin{equation}
    K = q^{H}\, , \quad K' = q^{H'} \, ,
\end{equation}	
	and 
    $ s_{q} (w) $
	is a 
    $q$-exponent
	which can be written in several forms
\begin{eqnarray}	
    s_{q}(w) & = & \prod_{n=0}^{\infty} ( 1 + q^{2n+1} w ) = \\
	     & = & 1 + \sum_{k=0}^{\infty} \frac{(-1)^{k} 
		q^{\frac{n(n-1)}{2}} w^{k}}{(q-q^{-1})\ldots(q^{k}-q^{-k})} = \\
	     & = & \exp \sum_{k=1}^{\infty}
			\frac{(-1)^{k} w^{k}}{k (q^{k}-q^{-k})} \: .
\end{eqnarray}	
	
	The term 
    $q$-exponent 
	is the most appropriate for the second form. In the third form
	there enters the 
    $q$-deformed dilogarithm, which was explored in particular in
\cite{FV2,FK}.	
	The title
    $q$-exponent 
	is strongly supported by the property, first found in
\cite{Schutz}:		
	let
    $ u $, $ v $
	be a Weyl pair
\begin{equation}	
	u v = q^{2} v u \, ,
\end{equation}		
	then
\begin{equation}	
	s_{q}(u) s_{q}(v) = s_{q}(u+v) \, .
\end{equation}		
	Let us note, that in
\cite{FV2}		
	it was found, that one more property holds
\begin{equation}	
	s_{q}(v) s_{q}(u) = s_{q}(u+v+ q^{-1}uv) = 
		s_{q}(u) s_{q}(q^{-1}uv) s_{q}(v)\, ,
\end{equation}			
	which can be called a "pentagon identity" and is a quantum
	deformation of the corresponding property of dilogarithm
\cite{FK}.		
	Now we see one more deficiency of the definition of the universal 
    $R$-matrix, 
	besides the necessity of using the
    $ \log K $.
	The function
    $ s_{q}(w) $
	behaves badly for
    $ q $
	lying on the unit circle
    $ |q| =1 $.
	Indeed, for example, in the third form of
    $ s_{q}(w) $	
	we see the small denominators. Thus the expression for the
	universal
    $R$-matrix 
	asks for some mending.

\section{The main idea}

	The problem of defining the
    $ \log $
	of operator appears already in a simpler example of Weyl pair
    $ u $, $ v $.	
	It is easy to realize the defining relation for
    $ u $, $ v $	
	via the Heisenberg pair
    $ P $, $ Q $	
	with relation
\begin{equation}	
	[ \, Q , P \, ] = 2 \pi i \hbar I \, .
\end{equation}			
	Indeed, the pair
\begin{equation}	
    u = e^{\alpha P} \, , \quad v = e^{\beta Q}
\end{equation}				
	satisfies Weyl relations 
\begin{equation}	
	u v = q^{2} v u 
\end{equation}		
	with
    $ \ln q = \frac{\pi \alpha \beta \hbar}{i} $.
	Thus the pair
    $ (P, Q) $
	defines
    $ (u, v) $.	
	However the inverse is not true, partly because the
    $ \log u $ or $ \log v $
	are badly defined. More subtle fact is that the pair
    $ P $, $ Q $	
	defines a second Weyl pair
\begin{equation}	
    \tilde{u} = e^{\tilde{\alpha} P} \, , 
	\quad \tilde{v} = e^{\tilde{\beta} Q}
\end{equation}				
	with a different phase
    $ \tilde{q} $,
    $ \ln \tilde{q} = \frac{\pi \tilde{\alpha} \tilde{\beta} \hbar}{i} $,
	which commutes with the first pair if
\begin{equation}	
    \alpha \tilde{\beta} = \hbar \, , \quad
    \tilde{\alpha} \beta = \hbar \, ,
\end{equation}				
	so that in particular	
\begin{equation}	
    \alpha \beta = \frac{\hbar^{2}}{\tilde{\alpha} \tilde{\beta}} \, .
\end{equation}				
	What is less trivial is the fact that together the commuting pairs
    $ (u, v) $	
	and
    $ (\tilde{u}, \tilde{v}) $		
	define naturally 
    $ P $ and $ Q $		
	for generic
    $ q $.		
	This fact was discussed explicitely in
\cite{FLMP},
	but of course could be traced to earlier literature, in particular
	to A.~Connes monograph on noncommutative geometry
\cite{Connes}
	and paper
\cite{Rieffel}.
	
	We see more coincisely, that the algebra
    $ \BB $
	generated by
    $ P $, $ Q $	
	is factored into the product
\begin{equation}	
     \BB = \AAA_{q} \otimes \AAA_{\tilde{q}}
\end{equation}					
	of commuting factors, generated by
    $ (u, v) $	
	and
    $ (\tilde{u}, \tilde{v}) $		
	correspondingly. I used the term "factor" with full algebraic
	meaning: indeed, neither
    $ \AAA_{q} $	
	nor
    $ \AAA_{\tilde{q}} $		
	have nontrivial center. However, I did not introduce any
    $*$-structure, 
	so the connection with v-Neumann theory (see e.g.
\cite{Connes}) 
	is incomplete. In particular, the last formula must contain some
	closure. In other words, the use of tensor product in this formula
	is somewhat loose. Indeed,
    $ \AAA_{q} $	
	and
    $ \AAA_{\tilde{q}} $		
	are factors of the type
    $ \TWO_{1} $		
	for generic
    $ q $	
	because they are infinite dimensional but allow for the trace
\begin{equation}	
     \tr \left( \sum a_{mn} u^{m} v^{n} \right) = a_{00} \, ,
\end{equation}						
	which is equal to 1 for the unit operator. On the other hand,
    $ \BB $	
	is a factor
    $ \ONE_{\infty} $
	as it can be realized as algebra of all operators in
    $ L_{2} (\RR) $.
	
	I hope that now the idea what to do in Quantum Group case is clear
	--- to define the
    $ \log K $	
	one is to extend the algebra
    $ \AAA_{q} $,	
	adding to it some dual generators. The definition of "dual" is
	clear for the Weyl type operators
\begin{equation}	
	\tilde{K} = (K)^{1/\tau}
\end{equation}							
	if we put 
    $ q = e^{i \pi \tau}$.	
	However not all generators of Quantum Group are of Weyl type. So we
	are to seek for the new set of generators which have this property.
	Fortunately, the theory of integrable models, which previously led
	to main relation of Quantum Group
\cite{KR},
	produces also this relevant set of generators. Indeed, the Quantum
	Group generators appeared first as the elements of the Lax operator
	of XXZ model. The lattice Sine-Gordon model introduced in 
\cite{IK}	
	belongs to this class (see e.g.
\cite{LesHouch}) 
	and in its turn naturally uses the Weyl-type generators. In the
	next section we shall give corresponding formulas in somewhat
	purified form.
	
\section{The explicite construction}
	Consider the algebra
    $ \CCC $	
	generated by four generators
    $ w_{1} $, $ w_{2} $, $ w_{3} $, $ w_{4} $,
	which is convenient to label by index
    $ n \in \ZZ_{4}$,
	so that
\begin{equation}	
	w_{n+4} = w_{n} \, .
\end{equation}							
	We impose Weyl-type relations on the nearest neighbours
\begin{equation}	
	w_{n} w_{n+1} = q^{2} w_{n+1} w_{n} \, ,
\end{equation}								
	whereas
\begin{equation}	
	w_{n} w_{m} = w_{m} w_{n} \, \quad |m-n| > 1 \, .
\end{equation}	
	Algebra
    $ \CCC_{q} $		
	(supplied by label
    $ q $) 
	has two central elements
\begin{equation}	
	Z_{1} = w_{1} w_{3} \, ; \quad Z_{2} = w_{2} w_{4} \, .
\end{equation}		
	It is a simple exercise to check that
\begin{eqnarray}
    e = i \frac{w_{1}+ w_{2}}{q-q^{-1}} \, , &&
    f = i \frac{w_{3}+ w_{4}}{q-q^{-1}} \, , \\
    K = q^{-1} w_{2} w_{3} \, , && 
    K' = q^{-1} w_{4} w_{1}
\end{eqnarray}	
	satisfy the defining relation of the Quantum Group. The relation
	between the central elements is as follows
\begin{equation}	
	J = Z_{1} Z_{2} \, ; \quad C = Z_{1} + Z_{2} 
\end{equation}			
	and so
    $ \CCC_{q} $		
	is a double cover of
    $ \AAA_{q} $.

	Following the reasoning of the previous section I introduce a
	second algebra
    $ \CCC_{\tilde{q}} $
	with generators
\begin{equation}	
	\tilde{w}_{n} = w_{n}^{1/\tau} \, .
\end{equation}			
	Both can be described in terms of the Heisenberg type generators
    $ p_{n} $	
	with relations\footnote{We put now the Planck constant
    $ \hbar $ equal to 1.}
\begin{equation}	
	[ \, p_{n} , p_{n+1} \, ] = - 2 \pi i I
\end{equation}				
	if we put
\begin{equation}	
    w_{n} = e^{bp_{n}} \, , \quad \tilde{w}_{n} = e^{p_{n}/b} \, ,
\end{equation}			
	where
\begin{equation}	
    q = e^{i \pi b^{2}} \, , \quad \tau = b^{2} \, .
\end{equation}	
	In particular, we see that
\begin{eqnarray}
    K = e^{b(p_{2}+p_{3})} \, , &&
    K' = e^{b(p_{1}+p_{4})} \, , \\
    \tilde{K} = e^{\frac{(p_{2}+p_{3})}{b}} \, , &&
    \tilde{K}' = e^{\frac{(p_{1}+p_{4})}{b}} \, ,
\end{eqnarray}		
	and the first factor in the universal 
    $R$-matrix 
	is expressed as
\begin{equation}	
    q^{ - \frac{H \otimes H'}{2}} = 
	e^{\frac{\pi}{2i} (p_{2} + p_{3}) \otimes (p_{1} + p_{4})} =
		\tilde{q}^{-\frac{\tilde{H} \otimes \tilde{H}'}{2}} 
\end{equation}	
	and does not depend on
    $ q $,
	serving both quantum groups
    $ \AAA_{q} $ and
    $ \AAA_{\tilde{q}} $.
	
	Let us turn now to the second factor in the universal 
    $R$-matrix. We have
\begin{equation}	
    s_{q} (- (q-q^{-1})^{2} e \otimes f) =
	s_{q} ( (w_{1} + w _{2} ) \otimes ( w_{3} + w_{4} ) ) =
\end{equation}		
	--- using Schutzenberger relation ---
\begin{equation}	
    = s_{q} ( w_{1} \otimes w_{3} ) s_{q} ( w_{1} \otimes w_{4} )
	s_{q} ( w_{2} \otimes w_{3} ) s_{q} ( w_{2} \otimes w_{4} ) \, ,
\end{equation}
	so that only Weyl-type combination enter here. Now I use the
	observation from
\cite{Varenna,FLMP}:
	consider the function
\begin{equation}
    \psi(p) = \exp \frac{1}{4} \oint_{-\infty}^{\infty}
	\frac{e^{ip\xi /\pi}}{\sh b\xi \sh \xi /b} \frac{d\xi}{\xi} \, ,
\end{equation}	
	where the singularity at
    $ \xi = 0 $
	in the integral is circled from above. It is easy to see that
\begin{equation}
    \psi(p) = \frac{s_{q}(w)}{s_{\tilde{q}}(\tilde{w})} \, ,
\end{equation}		
	where
    $ w = e^{bp} $,
    $ \tilde{w} = e^{p/b} $.
	Thus the integral
    $ \psi(p) $
	unites both
    $q$-exponents
   $ s_{q}(w) $ and $ s_{\tilde{q}}(\tilde{w}) $
	dual to each other. On the other hand it definitely has no their
	deficiencies, in particular, it does not suffer from the problem of
	the small denominators. The pentagon relation for it takes the
	form
\begin{equation}
    \psi(P) \psi(Q) = \psi(Q) \psi(P+Q) \psi(P)
\end{equation}			
	if
\begin{equation}
	[ \, P , Q \, ] = - 2 \pi i I \, .
\end{equation}			
	
	Function
    $ \psi(p) $	
	was used extensively in
\cite{FLMP,Varenna}
	and later by Kashaev in his new proposal for the knot invariants
\cite{Kknots}	
	and quantizing of the Techmuller space
\cite{KTech},
	which was done also independently by Chekhov and Fock
\cite{Fock}.
	
	The relation of
    $ \psi(p) $	
	to 
    $q$-exponents
	makes our proposal quite clear: consider the algebra
\begin{equation}	
    \DD = \CCC_{q} \otimes \CCC_{\tilde{q}} \, ,
\end{equation}
	generated by the generators
    $ w_{n} $, $ \tilde{w}_{n} $, $ n = 1, \ldots 4$.
	Quantum groups
    $ \UU_{q} $ and $ \UU_{\tilde{q}} $
	are naturally imbedded into
    $ \DD $
	by means of defining relations between Chevalley generators and
    $ w $-s.
	Algebra
    $ \DD $	
	can also be considered as being generated by the generators
    $ p_{n} $, $ n = 1, \ldots 4$.
	
	The element
    $ \RRR $
	in
    $ \DD \otimes \DD $
\begin{equation}	
    \RRR = \exp \left\{ \frac{\pi}{2i}(p_{2}+p_{3})\otimes(p_{1}+p_{4})\right\}
	\psi(p_{13}) \psi(p_{14}) \psi(p_{23}) \psi(p_{24}) \, ,
\end{equation}
	where
\begin{equation}	
    p_{ik} = p_{i} \otimes I + I \otimes p_{k}
\end{equation}
	plays the role of the universal
    $R$-matrix 
	for both
    $ \UU_{q} $ and $ \UU_{\tilde{q}} $.
	The Yang-Baxter relation
\begin{equation}	
    \RRR_{12} \, \RRR_{13} \, \RRR_{23} = \RRR_{23} \, \RRR_{13} \, \RRR_{12}
\end{equation}
	is an easy consequence of the pentagon relation, as was shown by
	R.~Kashaev and A.~Volkov
\cite{KV}.
	They also have obtained corresponding construction for the 
    $ SL_{q}(N) $
	case.
	
	There is a natural way to introduce the
    $*$-structure
	into
    $ \DD $:	
	one is to consider 
    $ p_{n} $
	to be selfadjoint
\begin{equation}	
	p_{n}^{*} = p_{n} \, .
\end{equation}	
	Corresponding formula in terms of 
    $ w $	
	looks as follows
\begin{equation}	
	w^{*} = w^{\overline{b}/b}\, ; \quad
	\tilde{w}^{*} = \tilde{w}^{b/\overline{b}}
\end{equation}	
	and does not respect the tensor structure of
    $ \DD $
	in general.

	However there are several particular cases when 
    $ * $-involution
	can be related to this structure:
	
	1.
    $ \tau > 0 $,
	so that 
    $ b $
	is real. The 
    $ * $
	takes the form
\begin{equation}	
	w^{*} = w \, ; \quad
	\tilde{w}^{*} = \tilde{w}
\end{equation}	
	and corresponds to 
    $ SL_{q}(2,\RR) $	
	reduction.
	
	2.
    $ \tau < 0 $,
	so that
    $ b $	
	is imaginary and 
\begin{equation}	
	w^{*} = w^{-1}\, ; \quad
	\tilde{w}^{*} = \tilde{w}^{-1}
\end{equation}		
	which corresponds to the
    $ SU_{q}(2) $	
	reduction.
	
	3.
    $ \tau = e^{i\theta} $,
	so that
    $ \tau = e^{i\theta/2} = \overline{\tilde{b}}$.	
	The involution takes form
\begin{equation}	
	w^{*} = \tilde{w}
\end{equation}			
	and so interchanges the factors in the modular double.
	
	For all three cases the parameter of the central extension of the
	class mapping group
\begin{equation}	
    C = 1 + 6 \left(\tau + 1 / \tau +2 \right) = 
	1 + 6 \left(b + 1/b \right)^{2}
\end{equation}
	is real. The values of 
    $ C $
	are
    $ C \geq 25$, $ C \leq 1 $
	and
    $ 1 \leq C \leq 25 $
	correspondingly for the case 1, 2 and 3. The relevance of this 
    $ * $
	operation to Conformal Field Theory is discussed in
\cite{FVK}.

	In the view of modular duality relation
    $ \tau \to -1/\tau $
	it is natural to call
    $ \DD $
	modular double.
	
	I conclude by proposing several problems:
	
	1. Give the generalization of our construction to other serieses of
	the Quantum Groups.
	
	2. Describe the coproduct in
    $ \UU_{q} $
	in terms of
    $w$-generators.
	
	3. Find a natural definition of the closure entering formal tensor
	products like those for algebras
    $ \BB $
	or
    $ \DD $.

	It is my pleasure to aknowledge the elucidating discussions with
	A.~Connes, A.~Volkov and A.~Vershik.

\end{document}